\documentclass[12pt]{amsart}

\usepackage{amsfonts,amsthm,amscd}

 \newtheorem{theorem}{Theorem}[section]

 \newtheorem{corollary}[theorem]{Corollary}

 \newtheorem{lemma}[theorem]{Lemma}

 \newtheorem{proposition}[theorem]{Proposition}

 \newtheorem{remark}[theorem]{Remark}

 \newtheorem{mtheorem}[theorem]{Main Theorem}

 \newtheorem{observation}[theorem]{Observation}

 \theoremstyle{definition}

 \newtheorem{defn}{Definition}

 \theoremstyle{remark}

 \newcommand{\nc}{\newcommand}

 \newcommand{\be}{\begin{equation}}

 \newcommand{\ee}{\end{equation}}

 \newcommand{\bea}{\begin{eqnarray}}

 \newcommand{\eea}{\end{eqnarray}}

 \newcommand{\bc}{\begin{center}}

 \newcommand{\ec}{\end{center}}

 \nc{\bmth}{\begin{mtheorem}} \nc{\emth}{\end{mtheorem}}

 \nc{\bth}{\begin{theorem}} \nc{\eth}{\end{theorem}}

 \nc{\bpr}{\begin{proposition}} \nc{\epr}{\end{proposition}}

 \nc{\ble}{\begin{lemma}}

  \nc{\ele}{\end{lemma}}

 \nc{\bco}{\begin{corollary}} \nc{\eco}{\end{corollary}}

 \nc{\bre}{\begin{remark}} \nc{\ere}{\end{remark}}

 \nc{\bob}{\begin{observation}} \nc{\eob}{\end{observation}}

\nc{\beq}{\begin{equation}}
\nc{\beqn}{\begin{equation*}}
\nc{\eeq}{\end{equation}}
\nc{\eeqn}{\end{equation*}}

 \nc{\f}{\frac} \nc{\rw}{\rightarrow} \nc{\To}{\longrightarrow}

 \nc{\Rw}{\Rightarrow}

 \nc{\nt}{\stackrel{\sim}{\nabla}} \nc{\De}{\Delta}

 \nc{\na}{\nabla} \nc{\al}{\alpha}

  \nc{\bet}{\beta}

  \nc{\va}{\vartheta}

 \nc{\ga}{\gamma} \nc{\G}{\Gamma}

  \nc{\la}{\lambda}

  \nc{\La}{\Lambda}

   \nc{\cS}{\mathcal{S}}

 \nc{\si}{\sigma} \nc{\de}{\delta} \nc{\ep}{\varepsilon}

 \nc{\ei}{\ep _i}

\nc{\eij}{\ep_{i,j}}

   \nc{\ea}{\ep _{\al}}

  \nc{\ej}{\ep _j}

  \nc{\eb}{\ep _{\bet}}

 \nc{\Z}{\mathbb{Z}} \nc{\Te}{\Theta}

 \nc{\om}{\omega}\nc{\Om}{\Omega} \nc{\ro}{\rho}

  \nc{\Or}{\mathcal{O}} \nc{\cA}{\mathcal{A}} \nc{\cD}{\mathcal{D}}

 \nc{\C}{\mathbb{C}^\star} \nc{\q}{\mathbf{Q}}

 \nc{\Q}{\mathbb{Q}}\nc{\bD}{\mathbb{D}} \nc{\bA}{\mathbb{A}}
\nc{\R}{\mathbb{R}}
 \nc{\Co}{\mathbb{C}} \nc{\cH}{\mathcal{H}}
 \nc{\cP}{\mathcal{P}}

 \nc{\cF}{\mathcal{F}}\nc{\cN}{\mathcal{N}}
\nc{\cf}{{\it \textsf{f}}}
 \nc{\ca}{\emph{\textbf{a}}} \nc{\bU}{\emph{\textbf{U}}}

\nc{\cM}{\mathcal{M}}\nc{\cL}{\mathcal{L}}\nc{\cC}{\mathcal{C}}
\nc{\cg}{\flat}
  
\nc{\ti}[1]{\tilde{#1}}
 
\nc{\Zpn}{\Z_p^n} \nc{\Cpn}{\Co_p^n}
\nc{\cZstar}{\mathcal{Z}^\times_p}
 \nc{\cX}{\mathcal{X}}  \nc{\cY}{\mathcal{Y}}
\nc{\cB}{\mathcal{B}}\nc{\Zpt}{\Z_p^\times}

 \nc{\fXY}{f_{(X,Y),\tau}}
\nc{\SL}{\mathrm{SL}_2}

\newcommand{\ra}{\rightarrow}

\newcommand{\ot}{\otimes}
\newcommand{\mtc}{\mathcal}

\newcommand{\Lam}{\Lambda}

\newcommand{\eps}{\epsilon}
\newcommand{\bn}{\begin}

\newcommand{\sub}{\subsection}

\newcommand{\D}{\Delta}

\numberwithin{equation}{section}
\newtheorem{thm}[equation]{Theorem}

\newcommand{\dw}{\downarrow}
\newcommand{\uw}{\uparrow}

\newcommand{\ch}{\chi}
\newcommand{\mtr}{\mathrm}

\title[Fusion partitions]{Fusion rings arising from normal Hopf subalgebras}
\begin{document}
\author{Sebastian Burciu and Vicentiu Pasol}
\address{Inst.\ of Math.\ ``Simion Stoilow" of the Romanian Academy
P.O. Box 1-764, RO-014700, Bucharest, Romania}
\date{\today}
\email{sebastian.burciu@imar.ro, vicentiu.pasol@imar.ro}


\begin{abstract}
For any normal commutative Hopf subalgebra $K=k^G$ of a semisimple Hopf algebra we describe the ring inside $kG$ obtained by the restriction of $H$-modules. If $G=\Z_p$ this ring determines a fusion ring and we give a complete description for it. The case $G=\Z_{p^n}$ and some other applications are presented.
\end{abstract}

\maketitle

\section{Introduction}

The character theory for normal Hopf subalgebras developed in \cite{cos} shows that the restriction functor from a semisimple Hopf algebra to a normal Hopf subalgebra has properties similar to that of a restriction from a group to a normal subgroup. Using this similarity we construct a $\Z_+$ -based ring attached to any normal commutative Hopf subalgebra. These rings are obtained as the image of the restriction map from the character ring of the larger Hopf algebra into the character ring of the normal Hopf subalgebra. If the normal Hopf subalgebra is of type $k^G$ where $G$ is a finite group group they determine special partitions on the group elements of $G$ that we called unital partitions. If $G$ is cyclic with a prime number $p$ of elements we complete classify these unital partitions and therefore their associated rings which in this case become fusion rings. We use our results to the study of certain fusion categories determined by the finite dimensional representations of a semisimple Hopf algebra. 

Semisimple Hopf algebras of certain small dimensions were studied extensively in
the literature and numerous classification results and non-trivial examples of Hopf algebras were obtained (see for example \cite{AN}, \cite{masf}, \cite{EG1},  \cite{G}, \cite{N2}, \cite{Nr}, \cite{Nar}, \cite{N5}, \cite{IK}). Recently group-theoretical fusion categories were introduced in \cite{ENO} and \cite{ENO2}. It was conjectured for several years that any Hopf algebra is group theoretical.  A class of Hopf algebras that are not group theoretical was described in \cite{N}. Their representation categories are described as equivariantization of certain group theoretical fusion categories. This brought up to the study a larger class of categories called weakly group theoretical \cite{ENO2}. The results from \cite{DGNO} and \cite{ENO2} enabled to classify up to a twist a larger class of Hopf algebras especially those of dimension $pq^2$ (see also \cite{JL}) and dimension $pqr$. Despite of that, a precise correspondence between the Hopf algebras classified up to isomorphism and the corresponding class of fusion categories is missing from the actual literature. We show that our methods are very useful in this direction. As an illustration, we completely describe the Grothendieck group of the unique semisimple Hopf algebra of dimension $2p^2$ with $p^2$ grouplike elements \cite{masf}. From here we easily obtain the Grothendieck group of the class of non group theoretical Hopf algebras of dimension $4p^2$ constructed in \cite{N}.

We work over an algebraically closed field $k$ of characteristic zero. For a vector space $V$ by $|V|$ is denoted the dimension of $V$. All the modules are supposed to be left modules. The comultiplication, counit and antipode of a Hopf algebra are denoted by $\Delta$, $\epsilon$ and $S$, respectively. We use Sweedler's notation $\D(x)=\sum x_1\ot x_2$ for all $x\in H$. All the other notations for Hopf algebras are those used in \cite{Montg}.
\section{Normal Hopf subalgebras}

Throughout this paper $H$ will be a semisimple (hence finite dimensional) Hopf algebra over an algebraically closed field $k$ of characteristic zero. Then $H$ is also cosemisimple and $S^2=\mtr{Id}$ \cite{Lard}.

The notation $\Lam_H \in H$ is used for the idempotent integral of $H$ ( $\eps(\Lam_H)=1$) and $t_H \in H^*$ for the idempotent integral of $H^*$ ($t_H(1)=1)$. Denote by $\mtr{Irr}(H)$ the set of irreducible characters of $H$ and let $C(H)$ be the character ring  of $H$. Then $C(H)$ is a semisimple subalgebra of $H^*$ \cite{Z} and $C(H)=\mtr{Cocom}(H^*)$, the space of cocommutative elements of $H^*$.


If $M$ is an $H$-module with character $\chi$ then $M^*$ is also an $H$-module with character $\chi^*=\chi \circ S$. This induces an involution $``\;^*\;":C(H)\ra C(H)$ on $C(H)$.

If $\mtr{Irr}(H)$ is the set of irreducible $H$-modules then from \cite{Montg} it follows that the regular character of $H$ is given by the formula
\begin{equation}
\label{f1}|H|t_{ _H}=\sum_{\ch \in \mtr{Irr}(H) }\chi(1)\chi.
\end{equation}

There is an associative nondegenerate bilinear form on $C(H)$ given by $$(\ch,\;\mu)=\ch\mu(\Lambda_H). $$ It follows that $(\ch,\;\mu)=m_{ _H}(\ch,\; \mu^*)$ where $m_{ _H}$ is the multiplicity form on $C(H)$. For two modules $M$ and $N$ of $H$ one has $m_{ _H}(\ch_M,\;\ch_N)=\mtr{dim}_k\mtr{Hom}_H(M,\;N)$. The pairs $\{\ch\}_{\ch \in \mtr{Irr}(H)}$ and $\{ \ch^*\}_{\ch \in \mtr{Irr}(H)}$ form dual bases with respect to $(\;,\;)$.

The following properties of $m_{ _H}$ will be used later:
$$m_{ _H}(\ch, \;\mu\nu)=m_{ _H}(\mu^*,\;\nu\ch^*)=m_{ _H}(\nu^*,\;\ch^*\mu)\;\;\text{and}$$ $$m_{ _H}(\ch,\mu)=m_{ _H}(\mu,\;\ch)=m_{ _H}(\mu^*,\;\ch^*)$$ for all $\ch,\mu,\nu \in C(H)$ (see Theorem 10 of \cite{NRi}).

This show that the $\Z$-module spanned by the irreducible characters is a fusion ring by the definition given below.

\sub{ Fusion rings}Let $\R_+$ be the semi-ring of non-negative real numbers.
We follow the definition from \cite{O} but allowing the structure constants to be any non-negative numbers. Let $R$ be a ring with identity which
is a finite rank $\Z$-module. A $\R_+$-basis of $R$ is a basis $B = \{X_i\}_{i=1,l}$ such that
$X_iX_j = \sum_{k=1}^l N^k_{ ij} X_k$, where $N^k_{ij} \in \R_+$. An element of $B$ will be called basic.
Define a non-degenerate symmetric $\Z$-valued inner product $m_R$ on $R$ as follows.
For all elements $X = \sum_{i=1}^l a_iX_i$ and $Y = \sum_{i=1}^l b_iY_i$ of $R$ we set
$$m_R(X, Y ) =\sum_{i=1}^l a_ib_i$$

\begin{defn} A  fusion ring or (a unital based ring) is a pair $(R, B)$ consisting of
a ring $R$ with a $\R_+$-basis $B = \{X_i\}_{i=1,l}$ satisfying the following properties:

(1) There exists $i_0 \in I$ such that $1 =  X_{i_0}$.

(2) There is an involution $i \mapsto i^*$ of $R$  such that the induced map $X = \sum_{i=1}^l a_iX_i \mapsto  X* =\sum_{i=1}^l a_iX_{i^*}$ satisfies
$$m_R(XY, Z) = m_R(X,ZY^ *) = m_R(Y, X^*Z)$$
for all $X, Y,Z \in R$.
\end{defn}

The second condition is equivalent to $N^k_{ij}=N^{i^*}_{jk^*}=N^{j^*}_{k^*i}$ for all $i,j, k$.

\subsection{Normal Hopf subalgebras}
Let $K$ be a Hopf subalgebra of $H$. Then $K$ is also semisimple and cosemisimple \cite{Montg}.
The Hopf subalgebra $K$ is called normal if $\sum h_1xSh_2 \in K$ for all $h \in H$ and $x \in K$. The restriction functor from $H$-modules to $K$-modules induces an algebra map $$res:C(H) \ra C(K).$$

Some results from \cite{cos} are recalled. On the set of irreducible $H$-characters there is an equivalence relation defined by $\ch \sim \mu$ if and only if their restriction to $K$ have a common constituent. This is the equivalence relation $r^{H^*}_{ _{L^*,\;k}}$ from \cite{cos} where $L=H//K$. It is proven that $\ch \sim \mu$ if and only if $\frac{\ch\dw_K^H}{\ch(1)}=\frac{\mu\dw_K^H}{\mu(1)}$ (see Theorem 4.3 of \cite{cos}). Thus the restriction of two irreducible $H$-characters $\ch$ and $\mu$ to $K$ either have the same irreducible constituents or they do not have common constituents at all.

The above equivalence relation determines an equivalence relation on the set of irreducible characters of $K$. For two irreducible $K$-characters $\al $ and $\beta$ one has $\al \sim \beta$  if and only if they are constituents of  $\ch\dw_{ _K}^H$ for some irreducible character $\ch$ of $H$. Let $\mtc{A}_0,\cdots \mtc{A}_{s}$ be the equivalence classes of the equivalence relation defined on $\mtr{Irr}(H)$. Let $\mtc{B}_0,\cdots \mtc{B}_{s}$ be the corresponding equivalence classes of the equivalence relation defined on $\mtr{Irr}(K)$.

For each $0 \leq  i \leq s$ put $b_i=\sum_{\al\in \mtc{B}_i}\al(1)\al \in C(K)$ and $a_i=\sum_{\ch \in \mtc{A}_i}\ch(1)\ch \in C(H)$.

The formula from subsection $4.1$ of \cite{cos} can be written as:

\bn{equation} \label{restrform}
\ch\dw_{ _K}^H=\frac{\ch(1)}{b_i(1)}b_i
\end{equation}

for all $\ch \in \mtc{A}_i$.
\bre\label{div}
Formula \ref{restrform} shows that $b_i(1)$ divides $\ch(1)$ in the case that all the characters in $\al\in\mtc{B}_i$ have degree one. 
\ere
On the other hand since the regular character of $H$ restricts to the regular character of $K$ multiplied by the index of $K$ in $H$ it follows that
\bn{equation} \label{restrform2}
a_i\dw_{ _K}^H=\frac{|H|}{|K|}b_i.
\end{equation}
In particular $a_i(1)=\frac{|H|}{|K|}b_i(1)$ for all $i$.

From the definition of the equivalence relation it follows that the set $\{\ch^*\;|\; \ch \in \mtc{A}_i\}$ is also a class of equivalence on $\mtr{Irr}(H)$. Therefore there is an involution $i\mapsto i^*$ on the index set of the equivalence classes such that  $\mtc{A}_{i*}=\{\ch^*\;|\;  \ch \in \mtc{A}_i\}$. Since restriction of characters commutes with taking the dual it follows that the same involution gives $\mtc{B}_{i*}=\{\al^*\;|\;  \al \in \mtc{B}_i\}$. Thus $a_i^*=a_{i^*}$ and $b_i^*=b_{i^*}$.

It was proven in \cite{cos} that the trivial character of $K$ forms by itself an equivalence class. Without loss of generality we will assume for the rest of the paper that $\mtc{B}_0=\eps_K$ and therefore $b_0=\eps_K$.

\begin{thm}\label{ident}Let $K$ be a normal Hopf subalgebra of $H$.
 With the above notations it follows that $$a_ia_j=\sum_{k=1}^sN^k_{ij}a_k$$ for some rational numbers $N^k_{ij}\geq 0$. Moreover $N^k_{ij}a_k(1)=N^{i^*}_{jk^*}a_i^*(1)$.
\end{thm}

\bn{proof}Let $L=H//K$ the quotient Hopf algebra and $t_L=\eps_K\uw^H_K$ the character of $H$ corresponding to the the trivial character of $K$ induced to $H$. Then $t_L$ is the integral of $L$. From Proposition \cite{cos} it follows that $a_i$ form a basis for the eigenspace of the operator of left multiplication with $t_L$ on $C(H)$ corresponding to the Frobenius-Perron eigenvalue $|L|$. Since $a_ia_j$ is also eigenvector for this operator it follows that $a_ia_j$ is a linear combination of the basis $a_k$. Comparing the multiplicities of irreducible characters in this product it follows that all the number $N^k_{ij}$ are rational. The second formula follows form the fact that $m_H(a_k,\; a_ia_j)=m_H(a_{i^*},\;a_ja_{k^*})$ and the obvious equality $m_H(a_i ,\;a_j)=\delta_{ij}a_i(1)$.
\end{proof}

\bco
 With the above notations it follows that $$b_ib_j=\sum_{k=1}^s\frac{N^k_{ij}}{l}b_k$$ where $l$ is the index of $K$ in $H$.
\eco

\subsection{$K=k^G$} Suppose now that $K$ is commutative, thus $K=k^G$ is for some finite group $G$. Since $\mtr{Irr}(k^G)=G$ it follows that the group $G$ is partitioned in subsets $G=\mtc{B}_0\coprod \mtc{B}_1\coprod\ldots\coprod \mtc{B}_s$ where $\mtc{B}_0=\{1\}$, the unit of the group $G$. It follows that in this situation $b_i=\sum_{x \in \mtc{B}_i}x$. Note that in this situation $b_i(1)=\mtr{card}(\mtc{B}_i).$ The product formula from the previous corollary implies that $\frac{N^k_{ij}}{l}$ are integers and therefore all the fusion coefficients $N^k_{ij}$ are integers in this situation.

\bth
Suppose that $G$ is a finite group and $K=k^G$ is a normal Hopf algebra of a semisimple Hopf algebra $H$. In the above settings, assume that there is $r$ a non-negative integer such that $\mtr{card}(\mtc{B}_i)=r+1$ for any $1\leq i \leq s$.  Let $R$ be the $\Z$-span of the elements $\frac{1}{\sqrt{r+1}}b_i$ with $1 \leq i \leq s$. Then $R$ is a fusion ring with the usual multiplication of $\Z G$.
\eth

\bn{proof}
Since $a_i(1)=lb_i(1)$ it follows that $a_i(1)=a_j(1)$ for all $1\leq i, j \leq s$. Then Corollary \ref{ident} implies that $ N^k_{ij}=N^{i^*}_{jk^*}$ for all nonzero $i, j, k$ One has to divide $b_i$ by ${\sqrt{r+1}}$ in order for the multiplicity form $m$ to be the one defined in the definition of a fusion ring.
\end{proof}

In the next section we will characterize these fusion rings for $G=\Z_p$. It will be shown that the above assumption on the cardinality of the sets $\mtc{B}_i$ is automatically satisfied. We will also consider the case $G=\mu_{p^n}$. We also consider the case $G=\mu_{p^n}$.
\section{Fusion partitions}
We start by defining the fusion partitions.

\begin{defn}
Let $G$ be a (commutative) group. A partition $G=I\coprod A_0\coprod\ldots\coprod A_s$ of $G$ is called a {\it unital partition} if the following conditions are satisfied:
\begin{enumerate}
\item $I=\{1\}$
\item Let $a_i:=\sum_{x\in A_i}x \in \Z[G]$. Then, the $\Z$-module generated by $\{a_i\}_i$ is a subalgebra of $\Z[G]$.
\end{enumerate}
A unital partition is called {\it fusion partition} if the following additional conditions are satisfied:
\begin{enumerate}
\item For all $A_i\in\cP$, there exists $i^*$ such that $A_{i^*}=\{x^{-1}\mid x\in A_i\}$
\item
$$
n_{i,j}^k\cdot |A_k|=n_{i,k^*}^{j^*}\cdot |A_{j^*}|,
$$
where $a_i\cdot a_j=\sum_k n_{i,j}^k a_k$. 
\end{enumerate}
\end{defn}

Note that if $k^G$ is a normal Hopf subalgebra of $A$ then as in the previous section restriction of simple $A$-modules to $k^G$ determines a fusion partition on $G$. The last identity from the definition follows from Proposition \ref{ident}.

We call a fusion partition trivial if any of its sets has a single element. Note that in the above settings one gets the trivial fusion partition if and only if $k^G$ is central in $A$.

\begin{defn}
Let $G=\mu_m$ be the cyclic group with $m$ elements and let $d|m$ be a divisor of $m$. For a partition $\cP=\{A_0,\ldots, A_s\}$ of $G$, let $\cP_d$ to be the subset of $\cP$ formed by the sets $A_i$ which contain a primitive $d^{th}$ root.
\end{defn}

We have the following important lemma:

\ble\label{power}
If $\cP$ is unital partition of $G=\mu_m$, then all the sets of $\cP_d$ have the same cardinal. Moreover, if $A\in\cP_d$, then all the sets in $\cP_d$ are of the form $A^n:=\{x^n\mid x\in A\}$ with $(n,m)=1$.
\ele

\begin{proof}
Let $A\in\cP_d$ such that $|A|$ is minimal. Since $(n,m)=1$ we have that $|A|=|A^n|$. Our assertion is therefore equivalent to the fact that $A^n\in\cP_d$ for all $(n,m)=1$.

Therefore it is enough to prove our lemma for the particular case when $n$ is a prime number not dividing  $m$.

In this case, we use the well known formula $(\sum_{x\in A} x)^n \equiv \sum_{x\in A}x^n$ mod $n$, which means that in $\Z/n\Z[G]$ we have $a^n=a(n)$, where $a=\sum_{x\in A}x$ and $a(n):=\sum_{y\in A^n}y$. On the other hand, our hypothesis shows that $a^n=\sum_i N_i a_i$. Reducing mod $n$ we get that $a(n)=\sum_i N_i a_i$ in $\Z/n\Z[G]$. Hence, there exists at least one $i$ such that $N_i\ne 0$ in $\Z/n\Z[G]$ (in fact there exists an $i$ such that $N_i= 1$ mod $n$. This means that in the equality  $a^n=\sum_i N_i a_i$ taken in $\Z/n\Z[G]$, in the left hand side we have exactly $|A|$ elements whose coefficients are $1$ and in the right hand side we have at least $|A_i|\ge |A|$ elements whose coefficients are $1$. This proves that in fact $|A_i|=|A|$ and there exists only one such $i$. Since every element in $A^n$ has to be in some $A_i$, this proves that in fact $A^n=A_i$.
\end{proof}

The fact that every $\cP_d$ is generated in the above sense by any single contained set gives us the following:

\bco
With the above notations, $\cP_d\bigcap \cP_{d'}\ne\emptyset$ is equivalent with $\cP_d=\cP_{d'}$.
\eco

\subsection{Unital partitions for $G=\mu_{p^n}$}$\;$\\

We consider in this subsection the case of unital partitions for $G=\mu_{p^n}$, the cyclic group with $p^n$ elements, where $p$ is an odd prime. In this case, $(\Z/p^n\Z)^\times$ is cyclic also. Let $\al$ be a natural number such that its class is a generator for $(\Z/p^n\Z)^\times$. Also, for the sake of short notations, we will change a little bit our previous notations:
\begin{enumerate}
\item{-} For $0\le k\le n$, we denote by $B_k:=\{x^{p^k}\mid x\;\; a\;\; generator\;\; for\;\; G\}$.
\item{-} We denote by $\cP_k:=\{A\in\cP\mid A\cap B_k\ne\emptyset\}$.
In the old notations this is $\cP_{p^{n-k}}$.
\item{-} We define the positive integer $u_k$ by the following:
   Let $A\in \cP_k$ and $y\in A\cap B_k$. Then $u_k$ is the smallest positive integer such that $y^{\al^{u_k}}\in A$. It is easy to see that this definition doesn't depend on the choice of $A$, $y$ or even $\al$ and always $u_k\le \phi(p^{n-k})$.
\end{enumerate}

We have the following lemmata:

\ble With the above notations, we have that
$$y^{\al^{tu_k}}\in A,$$
\ele

for  all integers $t$.
\begin{proof}
This is obvious from our previous lemma by induction. It is true for $t=0,1$. Assume that $y^{\al^{(t-1)u_k}}, y^{\al^{tu_k}}\in A$. Take $m\equiv \al^{u_k}$ mod $p^n$. This shows that $y^{\al^{tu_k}}\in A^m\bigcap A$ therefore $A=A^m$, so $y^{\al^{(t+1)u_k}}\in A$.
\end{proof}

\ble
We have $u_k\mid \phi(p^{n-k})$.
\ele

\begin{proof}
Indeed, let $u_k=v\cdot s$ with $v=(u_k,\phi(p^{n-k}))$ and $(s,\phi(p^{n-k}))=1$. It is clear that $v\le u_k$.

Let $w$ such that $s\cdot w=1$ mod $\phi(p^{n-k})$. On the other hand, $y^{\al^v}=y^{(\al^{sw})^v}.$ This proves that $y^{\al^v}=y^{\al^{wu}}\in A$. By the minimality of $u_k$, we get that either $v=0$, which is impossible, or $v\ge u_k$. Therefore, $u_k=v\mid \phi(p^{n-k})$.
\end{proof}

Let $\bet_k=\al^{u_k}$.
Fix a set $A \in \cP_k$ and an element $y \in A$. Put $A_{i,k}:=A^{\al^i}$ for all $i \geq 1$. Lemma \ref{power} implies $A_{i,k} \in \cP$.

\ble\label{b_i}
Let $B_{0,k}:=\{y,y^\bet_k,\ldots,y^{\bet_k^{r_k}}\}$, where $r_k+1=\f{\phi(p^{n-k})}{u_k}$. Also, let $B_{i,k}:=B_{0,k}^{\al^i}$.
Then, $B_{i,k}\subseteq A_{i,k}$ and $B_{i,k}\bigcap B_{j,k}\ne \emptyset$ if and only if $i\equiv j$ mod $u_k$.
\ele

\begin{proof}
The first part is clear by the definition of $A_{i,k}$ and $B_{i,k}$ and by the previous lemmata.
Now, let $y^{\al^{i+tu_k}}=y^{\al^{j+su_k}}\in B_{i,k}\bigcap B_{j,k}$. Then, $\al^{i+tu_k}\equiv \al^{j+su_k}$ mod $p^{n-k}$. This implies that $i-j+(t-s)u_k\equiv 0$ mod $\phi(p^{n-k})$. Since $u_k|\phi(p^{n-k})$ we get that $i-j \equiv 0$ mod $u_k$. The reverse assertion is obvious.
\end{proof}

\ble\label{nonint}
Let $\cP$ be a unital partition of $G=\mu_{p^n}$. With the above notations if $B_{j,k}\subset A_{i,k}$ then $j=i \;mod\;u_k$.
\ele

\begin{proof}
Assume that $B_{j, k}\subseteq A_{i, k}$, where $i<j<i+u_k$ (we can always rearrange the notation to have this condition). In this case, $B_{j-i,k}\subseteq A_{0,k}$ since $B_{j-i,k}=B_{u_k-i+j,k}=B_{j,k}^{\al^{u_k-i}}\subseteq A_{i,k}^{\al^{u_k-i}}=A_{0,k}$. We have therefore that $y^{\al^{j-i}}\in A_{0,k}$ and $0<j-i<u_k$. This contradicts the minimality of $u_k$.
\end{proof}

\bco\label{nonint2}
Let $\cP$ be a unital partition of $G=\mu_{p^n}$ and suppose that $\cP_k\bigcap \cP_{l}=\emptyset$ for $l \neq k$. Then $B_{i,k}= A_{i,k}$ for all $i$ and $k$.
\eco
\subsection{The case $G=\mu_p$} In this subsection we will describe the fusion partitions of $G=\mu_p$. Note that in this case $B_0=G\setminus \{1\}$ and $B_1=\{1\}$. We denote by $u$ the number $u_0$ corresponding to $\cP_0$. Let also $r:=r_0$ and $\beta:=\beta_0$.

We will prove the following result:

\bth\label{p} Let $G=\mu_p$ and let $\al$ be a generator for $(\Z/p\Z)^\times$. Fix $u|(p-1)$ and $x\in G$ a generator. Then $\cP=\{I,A_0,\ldots, A_u\}$  is a unital partition and any unital partition is obtained in this way. Here:
\begin{itemize}
\item{} $r+1=(p-1)/u$ is the order of $\al^u=:\bet$.
\item{} $A_i=\{x^{\al^i},x^{\al^i\cdot\bet},\ldots,x^{\al^i\cdot\bet^r}\}$, for all $i=\overline{0,u-1}$ and $I=\{1\}$.
\end{itemize}
Moreover, any unital partition is fusion partition.
\eth

\begin{proof}
The fact that any unital partition of $G=\mu_p$ has to be of the form described in the theorem follows from the previous corollary. Next we will show that any such partition is a fusion partition.
It is clear from the definition of $A_i$ that $A_i=A_{i+u}$, i.e. the index depends only modulo $u$.

Let's observe that in $\Z/p\Z$, we have that $-1=\al^{(p-1)/2}$ so by the definition of $A_{i^*}$ we have that $i^*=i+\f{p-1}{2}$ mod $u$. This proves the first fusion condition. Note that $A_i=A_{i^*}$ if and only if $u|\f{p-1}{2}$.

We compute $a_i\cdot a_j$:
$$
a_i\cdot a_j=(\sum_{m=0}^r x^{\al^i\bet^m})\cdot (\sum_{n=0}^r x^{\al^j\bet^n})=\sum_{m,n}|T^t_{i,j}|\cdot x^t,
$$
where $T^t_{i,j}:=\{(m,n)\in \Z/(r+1)\Z\times\Z/(r+1)\Z\mid \al^i\bet^m+\al^j\bet^n=t\}.$

Now, if $t\ne 0$, then $t=\al^k$ for some $k$, therefore in this case $|T^t_{i,j}|=|M^k_{i,j}|$, where: $M^k_{i,j}:=\{(m,n)\in \Z/(r+1)\Z\times\Z/(r+1)\Z\mid \al^i\bet^m+\al^j\bet^n=\al^k\}.$ Let $n^k_{i,j}:=|M^k_{i,j}|.$

To prove that the partitions described in the theorem are unital is enough to prove that $n^k_{i,j}=n^{k+u}_{i,j}.$ But this is obvious since there is an obvious bijection $ M^k_{i,j}\To M^{k+u}_{i,j}$ by $(m,n)\To (m+1,n+1)$

To prove the remaining fusion condition, we need to prove that $n_{i,j}^k=n_{k^*,j}^{i^*}$.
We have:
$$
M_{k^*,j}^{i^*}:=\{(m,n)\in \Z/(r+1)\Z\times\Z/(r+1)\Z\mid \al^{k^*}\bet^m+\al^j\bet^n=\al^{i^*}\}.
$$
The equation inside the set becomes: $-\al^{k}\bet^m+\al^j\bet^n=-\al^{i}$, i.e. $\al^{i}\bet^{-m}+\al^j\bet^{n-m}=\al^{k}$ which is  $\al^{i}\bet^{m'}+\al^j\bet^{n'}=\al^{k}$, where:
$$
{m'\choose n'}=\begin{pmatrix} -1&0\\-1&1\end{pmatrix}{m\choose n}
$$
The bijection of the two sets is immediate by the fact that the matrix involved is invertible.\end{proof}
\subsection{The case $G=\mu_{p^n}$ with some additional assumptions}

\begin{defn}
Let $G$ be a group and $\cP$ a partition of $G$. Then the set $A\in\cP$ is called singular if $A$ consists of a single element, i.e. $|A|=1$.
\end{defn}

In this section we will study another case of fusion partitions and that is the case $G=\mu_{p^n}$ such that $\cP$ has an additional condition:
$$
(\textbf{p})\;\;\;\;\;\;\;\;\;\;\;\;\; |A|={\rm is\;\; a\;\;} p-{\rm power}\; \forall A\in\cP.\;\;\;\;\;\;\;\;\;\;\;\;\;\;\;\;\;\;\;\;\;\;\;\;\;\;
$$

Suppose that $A$ is a Hopf algebra of dimension $p^n$ and $k^G$ with $G=\mu_{p^n}$ is a normal Hopf subalgebra of $A$. By \cite{MW} it follows that $A$ is of Frobenius type. Then Remark \ref{div} implies that the cardinality of any set of the partition $\cP$ is a power of $p$. Therefore $\cP$ satisfies the assumption $(p)$.
We have the following:
\ble
With the above notations, assuming $(\textbf{p})$, we have that $\cP_i\cap\cP_j=\emptyset$ for all $i\ne j$. That is $\cP=\coprod_k \cP_k$. Moreover, each $u_k$ is of the form $u_k=\phi(p^{b_k})$ for some $1\le b_k\le n-k$.
\ele

\begin{proof}
Suppose that $\cP_{k_1}=\cP_{k_2}=\ldots =\cP_{k_s}$. Let also $v:=|\cP_{k_i}|$ be the number of sets contained in them. Each set $A\in\cP_{k_i}$ contains exactly $\phi(p^{n-k_i})/v$ elements which are of the form $x^{p^{k_i}}$ with $x$ a primitive $p^n-th$ root. Therefore, each set $A$ contains exactly $\sum_{i =1}^s \phi(p^{n-k_i})/v$ elements. We may and will assume that $k_1<k_2<\ldots<k_s$. Condition $(\textbf{p})$ reads:
\beq
p^{n-k_s-1}\cdot \f{(p-1)\sum_{i=1}^s p^{k_s-k_i}}{v}=p^m
\eeq
The numerator of the fraction in the left hand side is relatively prime with $p$, therefore the fraction has to be a non-positive power of $p$. In particular $(p-1)\sum_{i=1}^s p^{k_s-k_i}\mid v$. On the other hand, $v\mid \phi(p^{n-k_i})=(p-1)p^{n-k_i-1}$ One easily gets a contradiction if $s>1$.\\
The second part is immediate.\end{proof}

Recall the sets $B_{i,k}$ and $A_{i,k}$ defined in Lemma \ref{b_i}.  Corollary \ref{nonint2} implies that under the assumption $(p)$ one has $B_{i,k}=A_{i,k}$ for all $i$ and $k$.
\ble
For any $i\in \Z/u_k\Z$ it follows that $$A_{i,k}=\{x^{p^k(\al^i+sp^{b_k})}\;| s=\overline{0,p^{n-k-b_k}-1}\}$$
\ele

\begin{proof}
It is enough to show that  the sets $\{\al^i+sp^{b_k}\;| s=\overline{0,p^{n-k-b_k}-1}\}$ and $\{{\al^i\cdot\bet_k^t}\mid t=\overline{0,p^{n-k-b_k}-1}\}$ coincide modulo $p^{n-k}$. Since $\bet_k^t=1(mod \;p^{b_k})$clearly the second set is a subset of the first one. Having the same number of elements the two sets must coincide.
\end{proof}

Using the previous lemma we will give a different indexation for the partition sets in $\cP_k$.  Note that $\al^i$ is coprime to $p^{b_k}$. Let $j\geq 1$ be any integer coprime with $p$. Write $B_{j,k}:=\{x^{p^k(j+sp^{b_k})}\;| s=\overline{0,p^{n-k-b_k}-1}\}$. Note that there is $i$ with $0 \leq i \leq u_k-1$ such that $\al^i=j(mod\; p^{b_k})$. Thus $B_{j,k}=A_{i,k}$ and any set $A_{i,k}$ can be written in this form.

Let $$b_{(j,\;k)}:=\sum_{z\in B_{j,k}}z=\sum_{s=0}^{r_k}x^{p^k(j+sp^{b_k})}$$

Let also $c_k:=r_k+1=p^{n-k-b_k}$ be the cardinality of a set in $\cP_k$.

\bre\label{rem}

1)Note that $x^{p^k(i+sp^{b_k})}=x^{p^k(i+tp^{b_k})}$  if and only if $s=t(mod\;c_k)$.

2)It is not hard to verify that $B_{i,\;k}=B_{j,\;k}$ if and only if $i=j(mod\;p^{b_k})$. This also follows from Lemma \ref{b_i}.
\ere

\ble
If the partition $\cP$ as above with property  $(\textbf{p})$ is fusion then $b_k+k$ is an increasing sequence, i.e. $b_{k+1}+k+1\ge b_k+k$ for all $k=\overline{0,n-1}$.
\ele

\begin{proof} The above inequality is true for $k=n-1$ since $b_{n-1}+n-1=1+n-1=n\le n=n+b_n$(we actually have equality). So we may assume that $0\le k\le n-1$. We may also assume that $b_k\ge 2$, otherwise the inequality is automatically satisfied. We look at the following product
$$b_{(1\;,k)}b_{(p-1,\;k)}=\sum_{s,t=0}^{r_k}x^{p^k(1+p^{b_k}s)}\cdot x^{p^k(p-1+p^{b_k}t)}=
\sum_{s,t=0}^{r_k}x^{p^{k+1}(1+p^{b_k-1}(s+t))}.$$ Remark that all the terms in the sums are in $B_{k+1}$ and therefore the product has to be a linear combination of terms of type $b_{i,k+1}$ for some appropriate values of $i$.

If $b_{k+1}<b_{k}-1$ then all the terms of the sum are in $B_{1, k+1}$. On the other hand since $1+p^{b_k-1}(s+t)=1+p^{b_{k+1}}(p^{b_k-1-b_{k+1}}(s+t))$ we see that not all the terms of $B_{1,k+1}$ can appear, since the coefficient of  $p^{b_{k+1}}$ in the above expression is divisible by $p$. Thus $b_{k+1}\geq b_{k}-1$.
\end{proof}

\ble
Assume that  $\cP$ satisfies the conditions of the previous lemma. Then:
\bn{enumerate}
\item If $k<l$ then $b_{(i,\;k)}b_{(j,\;l)}=c_lb_{(i+p^{k-l}j,\;k)}$.
\item If $i+j=p^sm$ with $0\leq s <b_k$ and $(m,p)=1$ then $$b_{(i,\;k)}b_{(j,\;k)}=c_k\sum_{t=0}^{p^{b_{k+s}-b_k-s}-1}b_{(m+p^{b_k-s}t, \;k+s)}$$%
\item
    If $i+j=p^{b_k}m$  then $B_{i,k}^*=B_{j,k}$ and
    $$b_{(i,\;k)}b_{(j,\;k)}=c_k\sum_{b_k+k \leq l \leq n}\;\;\sum_{1\leq v \leq p^{b_l},\; (v,p)=1}b_{(v,\;l)}$$
\end{enumerate}
\ele
\begin{proof}

1) If $k<l$ then $$b_{(i,\;k)}b_{(j,\;l)}=\sum_{s=0}^{r_k}x^{p^k(j+p^{b_k}t)} \cdot \sum_{t=0}^{r_l}x^{p^l(j+p^{b_l}t)}=\sum_{s=0}^{r_k}\sum_{t=0}^{r_l}
    x^{p^k(i+p^{k-l}j+p^{b_k}(s+p^{l-k+b_l-b_k}t))}$$
    Note that $b_l+l\geq b_k+k$ and therefore $l-k+b_l-b_k\geq0$. It can be checked that the elements $s+p^{l-k+b_l-b_k}t$ cover all the elements of $\Z_{c_k}$ (each of them $c_l$ times) when $s$ and $t$ run as in the sum.
%

2) If $i+j=p^sm$ then  $$b_{(i,\;k)}b_{(j,\;k)}=\sum_{t, t'=0}^{r_k}x^{p^{k+s}(m+p^{b_k-s}(t+t'))}$$

First note that $m+p^{b_k-s}t$ is prime with $p$ since $s<b_k$. Each element of the set $B_{m+p^{b_k-s}t, k+s}$ appears in the above product. Indeed any element of this set is of the form $y=x^{p^{k+s}(m+p^{b_k-s}t+p^{b_{k+s}}u)}$ for some integer $u$ and then it follows that $y=x^{p^{k+s}(m+p^{b_{k}-s}(t+p^{b_{k+s}-(b_k-s)}u))}\in B_{i,k}B_{j,k}$.

Next we show that each element of $G$ that appears in the product $b_{(i,\;k)}b_{(j,\;k)}$ has multiplicity $c_k$. Indeed any such element is of the form $y=x^{p^{k+s}(m+p^{b_k-s}(t+t'))}$ for some $0\leq t, t' \leq c_k-1$. Item 1 of Remark \ref{rem} implies that $y$ depends on the class of $t+t'$ modulo $c_k$ and therefore it appears $c_k$ times.

In order to finish the proof it is enough to check that the equality $B_{m+p^{b_k-s}t, k+s}= B_{m+p^{b_k-s}t', k+s}$ with $t$ and $t'$ elements of $\Z_{c_k}$ holds if and only if $t=t'(mod\;p^{b_{k+s}-b_k+s})$. Note that $p^{b_{k+s}-b_k+s}|c_k$ since $b_{k+s}-b_k+s \leq n-k-b_k$.

By the second item of Remark \ref{rem} one has that $B_{m+p^{b_k-s}t, k+s}= B_{m+p^{b_k-s}t', k+s}$ if and only if $p^{b_k-s}t=p^{b_k-s}t'(mod\;p^{b_{k+s}})$ which is the same as $t=t'(mod\;p^{b_{k+s}-b_k+s})$.

3) If $i+j=p^{b_k}m$ then $$b_{(i,\;k)}b_{(j,\;k)}=\sum_{t, t'=0}^{r_k}x^{p^{k+b_k}(m+t+t')}.$$ It is clear that $m+t+t'$ covers each residue modulo $p^{n-k-b_k}$ exactly $c_k$ times as $t$ and $t'$ runs through all the residues modulo $c_k=p^{n-k-b_k}$. Since the identity element $1 \in G$ appears in the above product it follows that $B_{i,k}^*=B_{j,k}$.
\end{proof}

We prove now the following theorem which is a characterization for the fusion partition with property $(\textbf{p})$:

\bth
Fix $n\ge 1$ a natural number and let  $\{b_k\}_k$ be a sequence of numbers which satisfy the conditions:
\begin{enumerate}
\item $b_n=0$ and $b_{n-1}=1$
\item $b_{k+1} +1\ge b_k$, i.e. the finite sequence $b_k+k$ is increasing.
\end{enumerate}
Fix $x$ a primitive $p^{n}-th$ root of unity and $\al$ a natural number such that $\al$ mod $p^n$ is a generator of $(\Z/p^n\Z)^\times$. Put $u_k:=\phi(p^{b_k})$ and $\bet_k:=\al^{u_k}$.
Let $A_{i,k}:=\{x^{p^k\cdot\al^i\cdot\bet_k^t}\mid t=\overline{0,p^{n-k-b_k}}\}$ for all $k=\overline{0,n}$ and for all $i\in \Z/u_k\Z$.

Put $\cP_k:=\{ A_{i,k}\mid i\in \Z/u_k\Z\}$ and put $\cP:=\coprod_k \cP_k$. Then
$\cP$ is a fusion partition for $G=\mu_{p^n}$ which satisfies property $(\textbf{p})$ and every fusion partition of $G=\mu_{p^n}$ which satisfies property $(\textbf{p})$ is obtained in this way.
\eth

\begin{proof}
By the above lemmata, we know that every fusion partition of $G=\mu_{p^n}$ which satisfies property $(\textbf{p})$ has to be of the form described in the theorem. Note that $a_{i,k}\cdot a_{j,l}$ is a linear combination of $a_{q,m}$ by the previous Lemma which shows that the partition is unital. By the same lemma is straightforward to verify that $N^k_{ij}c_k=N^{i^*}_{jk^*}c_{i^*}=N^{j^*}_{k^*i}c_{j^*}$ for all set indices $i,j, k$.\end{proof}
\section{Non-group theoretical Hopf algebras of dimension $4p^2$}
Consider $H_p$ the non group theoretical Hopf algebra of dimension $4p^2$ constructed in \cite{N}. Using the above methods we will compute the Grothendieck group of $H_p$. In order to do this  we need to study fusion partition for $G=\mu_p\times\mu_p$ which have $p$ singularities and each nonsingular set has two elements.
\subsection{Some fusion partitions on $G=\mu_p\times\mu_p$} We will prove the following:

\bth
In the case mentioned above, there exist $\al,\bet\in G$ such that each singular set is of the form $A_i:=\{\al^i\}$ and each nonsingular set is of the form $B_{i,j}:=\{\al^i\bet^j,\al^i\bet^{-j}\}$ with $i\in\overline{0,p-1}$ and $j\in\overline{1,p-1}$.
\eth
\begin{proof}
The elements of $\cP_{sing}$ form a group with $p$ elements which we will denote by $G_0$. Let $\al \in G$ be its generator. If $A=\{\beta, \beta'\}$ is a set of cardinality two from the fusion partition then it can be shown as above $A^m:=\{\bet^m, \;\bet'^m\}$ is also a set of the given fusion partition. Suppose $\bet'=\gamma\beta^i$ for some $\gamma \in G_0$ and $1 \leq i \leq p-1$. It follows that $\gamma A^i:=\{\gamma\bet^i,\;\gamma^{i+1}\bet^{i^2}\}$ is also a set of the fusion partition. Therefore $\gamma A^i=A$ which gives that $\bet=\gamma^{i+1}\bet^{i^2}$. This implies $i^2=1(mod \;p)$. If $i=1 (mod \;p)$ then $\beta=\gamma^2\bet$ which gives $\gamma=1$ since $p$ is odd. But then $A=\{\bet,\;\bet\}$ which is impossible. It follows that $i=-1(mod \;p)$ and $A=\{\bet, \; \gamma\bet^{-1}\}$. Suppose $\gamma \neq 1$, thus $\gamma=\al^s$ for some $1 \leq s\leq p-1$. Choose $i$ such that $2i=-s (mod\; p)$. Then $\la^iA=\{\al^i\beta,\;\al^{i+s}\beta^{-1}\} $ and the inverse of the first element of the set is the other element of the set. Replacing $\beta$ by $\al^i\beta$ the theorem follows.\end{proof}
\subsection{Non-group theoretical Hopf algebras of dimension $4p^2$}
 From \cite{N} $H_p$ is an extensions of $A_p$ by $k\Z_2$ where $A_p$ is the unique semisimple Hopf algebra of dimension $2p^2$ with exactly $p^2$ group-like elements (see\cite{masf}). Thus $H_p$ fits into a short exact sequence of Hopf algebras
\beq
k \ra k\Z_2 \ra H_p \ra A_p \ra k.
\eeq
 As algebras one has that $
H_p \cong H^*_p \cong
k^{2p}\oplus M_2(k)^{\frac{p(p-1)}{2}}\oplus M_p(k)^2.$
\subsection{The Grothendieck group of $A_p$} First we compute the \\Grothendieck group of $A_p$. By \cite{masf} there is a
short exact sequence of Hopf algebras
\beq
k \ra k^{\Z_p\times \Z_p} \ra A_p \ra k\Z_2 \ra k.
\eeq
From the same paper it follows that $A_p \cong k^{2p} \oplus M_2(k)^{\frac{p(p-1)}{2} }$ and $A^*\cong k^{p^2} \oplus M_p(k)$ as algebras. Let $\psi$ be a generator of $G(A^*)\cong \Z_{2p}$. The subalgebra $K=k^{\Z_p\times \Z_p}$ is normal in $A_p$ since it is of index $2$ \cite{Na3}. Therefore it determines a fusion partition on $G=\Z_p\times \Z_p$. Since the length of each subset of the fusion partition divides the degree of an irreducible character of $A_p$ it follows that this length can be at most two. The trivial character induced from $K$ to $A$ has degree two so it is the sum of the trivial character and other one dimensional character of $A$.  By Corollary 2.5 from \cite{Bker} the irreducible constituents of  $\eps_{K}\uw^{A_p}$ form a fusion subcategory of $\mtr{Rep}(A_p)$. Thus the other one dimensional character has order two in $G(A_p^*)$  and therefore it should be $\psi^p$. It follows that the fusion partition on $\Z_p\times \Z_p$ has singularities (since not all one dimensional representations) restrict to the trivial representation. Clearly $K$ is not central in $A_p$. Therefore the fusion partition has exactly $p$ singularities and the results from the previous subsection can be applied.

Let $\al$ be the character generating the group of the singularities of the above fusion partition and $\beta$ another character of $K$. One may assume that $\psi\dw_K=\al$. Since $\psi^{p+r}\dw_K=\psi^{r}\dw_K$ one has that $\psi^{p+r}\ch=\psi^r\ch$ for any irreducible character $\ch$ of degree $2$. The description of the fusion partitions from the previous subsection implies that for any $1 \leq s \leq \frac{p-1}{2}$ there is a two dimensional characters $\ch_s$ of $A_p$ such that $\ch_s\dw_K=\beta^s+\beta^{-s}$. Since the restriction map is an algebra map it follows that

$$\ch_s\ch_t\dw_K=\eps+ \beta^{s+t}+\beta^{-(s+t)}+\beta^{s-t}+\beta^{-(s-t)}.$$

Looking at possible constituents of the product $\ch_s\ch_t$ with the above restriction to $K$ it follows that $\ch_s\ch_t=\ch_{s+t}+\ch_{s-t}$ for $s \neq t$ and $\ch_s^2=\eps+\psi^p +\ch_{2s}$ (all indices are denoted modulo $p$). This completely determines the Grothendieck group of $A_p$ since any $2$-dimensional character of $A_p$ is of the type $\psi^{l}\ch_s$ for some integers $0 \leq l \leq p-1$ and $1 \leq s \leq \frac{p-1}{2}$.

\subsection{The Grothendieck group of $H_p$} Since $A_p$ is a quotient Hopf algebra of $H_p$ it follows that $\mtr{Rep}(A_p) \subset \mtr{Rep}(H_p)$. $H_p$ has in addition $2$ more characters of degree $p$. Denote them by $ \eta_1$ and $\eta_2$. From the description of $Rep(H_p)$ as a $\Z_2$-equivariantization of a Tambara-Yamagami category \cite{N} it follows that $\eta_1^*=\eta_2$ and $\eta_2^*=\eta_1$. Let $L(\eta_i)$ the subgroup of $G(H_p^*)\cong \Z_{2p}$ consisting of the one dimensional representations of $H_p$ that appear in $\eta_1\eta_1^*$. Since other characters that can appear in this product have even degree it follows that  $L(\eta_i)$ has odd degree. Therefore $L(\eta_i)=<\psi^2>$ and $\psi\eta_1=\eta_2$, $\psi\eta_2=\eta_1$. A similar argument implies that $\eta_1\psi=\eta_2$, $\eta_2\psi=\eta_1$. By duality it follows from \cite{N} that
\beq
k \ra A_p^*  \ra H_p \ra k\Z_2 \ra k.
\eeq
is an extension of Hopf algebras. The category $\mtr{Rep}(A_p^*)$ is a Tambara-Yamagami category corresponding to group $H=\Z_p \times \Z_p$ and with a $p$ dimensional representation whose character will be denoted by $d$. Since $\eps_{A_p^*}\uw^{H_p}_{A_p^*}$ has degree two it is a sum of $\eps_{H_p}$ and another one dimensional representation. By Corollary 2.5 from \cite{Bker} the irreducible constituents of  $\eps_{A_p^*}\uw^{H_p}$ form a fusion subcategory of $\mtr{Rep}(H_p)$. Thus the other one dimensional representation has a character of order two in $G(H_p^*)$  and therefore it should be $\psi^p$. Let $g=\psi\dw_{A_p^*}$. Since $\ch_1\in \mtr{Rep}H_p$ is self dual it follows that $\ch_1\dw_{A_p^*}=h+h^{-1}$ for some $h \in H$. Then $\ch_s\dw_{A_p^*}=h^{s}+h^{-s}$ for all $1 \leq s \leq \frac{p(p-1)}{2}$. Thus $\psi^{i}\ch_s\dw_{A_p^*}=g^ih^{s}+g^ih^{-s}$. Since the regular character of $H_p$ restricts to twice the regular character of  $A_p^*$ it follows that $\eta_1\dw_{A_p^*}=\eta_2\dw_{A_p^*}=d$. By Frobenius reciprocity this implies that $d\uw_{A_p^*}^{H_p}=\eta_1+\eta_2$.
From Proposition 2 of \cite{Bd} it follows $\ch\mu\dw_{A_p^*}=(\ch\dw_{A_p^*}\mu)\uw^{H_p}$ for all $\ch, \mu$ characters of $H_p$. Applying the above formula for $\ch=\ch_s$ and $\mu=\eta_1$ it follows that $\ch_s\eta_1=2(\eta_1+\eta_2)$. Similarly $\ch_s\eta_2=\eta_1+\eta_2$. There are three possibilities that we list next:

1)  $\ch_s\eta_1=\ch_s\eta_2=\eta_1+\eta_2$,

2) $\ch_s\eta_i=2\eta_i$ for $i=1,2$,

3)  $\ch_s\eta_1=2\eta_2$ and $\ch_s\eta_2=2\eta_1$.

The last two possibilities are excluded by the following argument: in both situations $\ch_s^2\eta_1=4\eta_1$ and the formula for $\ch_s^2$ from the previous subsection implies that $\psi^{p} \in L(\eta_1)$ which is impossible since $L(\eta_i)=<\psi^2>$.

In order to complete the description of the Grothendieck group of $H_p$ one needs to describe $\eta_1\eta_2$ and $\eta_1^2$ and $\eta_2^2$. Since $\psi^i\ch_s\eta_i=\eta_1+\eta_2$ it follows that $m(\psi^i\ch_s, \eta_1\eta_2)=m(\eta_1,\;\psi^i\ch_s\eta_1)=1$. Thus each irreducible character of degree $2$ of $H_p$ appears with multiplicity $1$ in the product $\eta_1\eta_2$. This gives the formula

$$\eta_1\eta_2=\sum_{j=0}^{p-1}\psi^{2j}+\sum_{i=0}^{p-1}\sum_{s=1}^{p-1}\psi^i\ch_s.$$

Using the same formula from \cite{Bd} it follows that $\eta_1(\eta_1+\eta_2)=\eta_1d\uw^{H_p}=\sum_{i,j=0}^{p-1}g^ih^j\uw^{H_p}=\sum_{i=0}^{2p-1}\psi^{i}
+2\sum_{i=0}^{p-1}\sum_{s=1}^{p-1}\psi^i\ch^s.$ This implies that $$\eta_1^2=\sum_{j=0}^{p-1}\psi^{2j+1}+\sum_{i=0}^{p-1}\sum_{s=1}^{p-1}\psi^i\ch_s.$$
Similarly one gets the same formula for $\eta_2^2$.
This completes the description of the Grothendieck group of $H_p$.

\bre
\ere
The above computations show that the Grothendieck group of $H_p$ is commutative. Moreover, completely similarly one can obtain the the Grothendieck group of $H_p^*$. Thus $H_p$ and $H_p^*$ have the same Grothendieck group which implies that one is a twist of the other \cite{Nics}.
\bibliographystyle{amsplain}
\bibliography{fusion}

\providecommand{\bysame}{\leavevmode\hbox to3em{\hrulefill}\thinspace}
\providecommand{\MR}{\relax\ifhmode\unskip\space\fi MR }
\providecommand{\MRhref}[2]{%
  \href{http://www.ams.org/mathscinet-getitem?mr=#1}{#2}
}
\providecommand{\href}[2]{#2}
\begin{thebibliography}{10}

\bibitem{Bker}
S.~Burciu, \emph{\textnormal{Normal Hopf subalgebras of semisimple Hopf
  Algebras}}, to appear Proc.\ A. M. S.\, math.RA/0611068 (2006).

\bibitem{Bd}
\bysame, \emph{\textnormal{On some representations of the Drinfel'd double}},
  J. Algebra \textbf{296} (2006), 480–504.

\bibitem{cos}
\bysame, \emph{\textnormal{Coset Decomposition For Semisimple Hopf Algebras}},
  to appear, Comm. Alg (arXiv:0712.1719) (2007).

\bibitem{DGNO}
V.~Drinfeld, S.~Gelaki, D.~Nikshych, and V.~Ostrik,
  \emph{\textnormal{Group-theoretical properties of nilpotent modular
  categories,}}, arXiv:0704.0195v2.

\bibitem{EG1}
P.~Etingof and S.~Gelaki, \emph{\textnormal{Semisimple Hopf algebras of
  dimension pq are trivial}}, J. Algebra \textbf{210} (1998), no.~2, 664–669.

\bibitem{ENO2}
P.~Etingof, D.~Nikshych, and V.~Ostrik, \emph{\textnormal{Weakly
  group-theoretical and solvable fusion categories}}, arXiv:0809.3031.

\bibitem{ENO}
P.~Etingof, D.~Nikshych, and V.~Ostrik, \emph{\textnormal{On fusion
  categories}}, Annals of Mathematics \textbf{162} (2005), 581--642.

\bibitem{G}
S.~Gelaki, \emph{\textnormal{Quantum groups of dimension $pq^2$}}, Israel J.
  Math. \textbf{102} (1997), 227–267.

\bibitem{JL}
David Jordan and Eric Larson, \emph{\textnormal{On the classification of
  certain fusion categories}}, arXiv:0812.1603 (2008).

\bibitem{Lard}
R.~G. Larson and D.~E. Radford, \emph{\textnormal{Finite dimensional
  cosemisimple \textnormal{Hopf} Algebras in characteristic zero are
  semisimple}}, J. Algebra \textbf{117} (1988), 267--289.

\bibitem{IK}
H.~Kosaki M.~Izumi, \emph{\textnormal{Kac algebras arising from composition of
  subfactors: general theory and classification}}, Mem. Amer. Math. Soc.
  \textbf{158} (2002), no.~750.

\bibitem{masf}
A.~Masuoka, \emph{\textnormal{Some further classification results on semisimple
  Hopf algebras}}, Comm. Algebra \textbf{178} (1996), no.~21, 307--329.

\bibitem{Montg}
S.~Montgomery, \emph{\textnormal{Hopf algebras and their actions on rings}},
  vol.~82, 2nd revised printing, Reg. Conf. Ser. Math, Am. Math. Soc,
  Providence, 1997.

\bibitem{MW}
S.~Montgomery and S.~Witherspoon, \emph{\textnormal{Crossed products for Hopf
  Algebras}}, J. Pure and Appl. Algebra \textbf{111} (1988), 381--385.

\bibitem{AN}
S.~Natale N.~Andruskiewitsch, \emph{\textnormal{Examples of self-dual Hopf
  algebras}}, J. Math. Sci. Univ. Tokyo \textbf{6} (1999), no.~1, 181–215.

\bibitem{N5}
S.~Natale, \emph{\textnormal{Hopf algebra extensions of group algebras and
  Tambara-Yamagami categories}}, arXiv:0805.3172.

\bibitem{Na3}
\bysame, \emph{\textnormal{On semisimple Hopf algebras of dimension $pq^2$}},
  J. Algebra \textbf{221} (1999), no.~1, 242--278.

\bibitem{N2}
\bysame, \emph{\textnormal{On semisimple Hopf algebras of dimension $pq^2$,
  II}}, Algebr. Represent. Theory \textbf{3} (2001), no.~3, 277--291.

\bibitem{Nar}
\bysame, \emph{\textnormal{On semisimple Hopf algebras of dimension $pq^r$}},
  Algebr. Represent. Theory \textbf{7} (2004), no.~2, 173--188.

\bibitem{Nr}
\bysame, \emph{\textnormal{Semisolvability of semisimple Hopf algebras of low
  dimension}}, no. 186, Mem. Am. Math. Soc., Am. Math. Soc., Providence, RI,
  2007.

\bibitem{NRi}
W.~D. Nichols and M.~B. Richmond, \emph{\textnormal{The Grothendieck group of a
  Hopf algebra}}, Journal of Pure and Applied Algebra \textbf{106} (1996),
  297--306.

\bibitem{Nics}
D.~Nikshych, \emph{\textnormal{$K_0$-rings and twistings of finite dimensional
  semisimple Hopf algebras}}, Comm. Algebra \textbf{26} (1998), 321--342.

\bibitem{N}
D.~Nikshych, \emph{\textnormal{Non group-theoretical semisimple Hopf algebras
  from group actions on fusion categories}}, Selecta Math. \textbf{14} (2009),
  no.~1, 145 -- 161.

\bibitem{O}
V~Ostrik, \emph{\textnormal{Module categories, weak Hopf algebras and modular
  invariants}}, Transform. Groups \textbf{26} (2003), no.~8, 177--206.

\bibitem{Z}
Y.~Zhu, \emph{\textnormal{Hopf algebras of prime dimension}}, Int. Math. Res.
  Not. \textbf{1} (1994), 53--59.

\end{thebibliography}
\end{document}